\documentclass[11pt]{article}

\usepackage{graphicx}
\usepackage{amssymb, amsmath,}
                                                                                                                                                                                                                                                                                                                                                                                                                                                                                                                                                                                                                                                                                                                                                                                                                                                                                                                                                                                                                                                                                                                                               \newtheorem{theorem}{Theorem}[section]

\newtheorem{definition}[theorem]{Definition}
\newtheorem{lemma}[theorem]{Lemma}
\newtheorem{Remark}[theorem]{Remark}

\newtheorem{Conjecture}[theorem]{Conjecture}
\title{\textbf{Intersection of two quadrics with no common hyperplane in $\mathbb{P}^{n}(\mathbb{F}_q)$}}

\date{}
\author{$\text{Fr\'ed\'eric A. B. Edoukou}^{*}$, $\text{San Ling}^{*}$, Chaoping Xing{\thanks{This work is supported by
MOE-AcRF Tier 2 Research Grant, Singapore (No.T206B2204).}}\\
Division of Mathematical Sciences,\\  Nanyang Technological University, \\
  21 Nanyang Link,
 Singapore 637371. \\
 E.mail : $\{\mathrm{abfedoukou, lingsan, xingcp} \}$@ntu.edu.sg}
\begin{document}
\maketitle

{\footnotesize \begin{flushleft} \textbf{Abstract}
\end{flushleft}
Let $\mathcal{Q}_1$ and $\mathcal{Q}_2$ be two arbitrary quadrics
with no common hyperplane in ${\mathbb{P}}^n(\mathbb{F}_q)$. We give
the best upper bound for the number of points in the intersection of
these two quadrics. Our result states that $| \mathcal{Q}_1\cap
\mathcal{Q}_2|\le 4q^{n-2}+\pi_{n-3}$. This result inspires us to
establish the conjecture on the number of points of an algebraic set
$X\subset {\mathbb{P}}^n(\mathbb{F}_q)$ of dimension $s$
 and degree $d$: $|X(\mathbb{F}_q)|\le dq^s+\pi_{s-1}$.  
}
%%%%%%%%%%%%%%%%%%%%%%%%%%%%%%%%%%%%%%%%%%%%%%%%%%%%%%%%%%%%%%%%%%%%%%%%%%%%%%%%%%%%%%%
%%%%%%%%%%%%%%%%%%%%%%%%%%%%%%%%%%%%%%%%%%%%%%%%%%%%%%%%%%%%%%%%%%%%%%%%%%%%%%%%%%%%%%%%
%%%%%%%%%%%%%%%%%%%%%%%%%%%%%%%%%%%%%%%%%%%%%%%%%%%%%%%%%%%%%%%%%%%%%%%%%%%%%%%%%%%%%%%
%%%%%%%%%%%%%%%%%%%%%%%%%%%%%%%%%%%%%%%%%%%%%%%%%%%%%%%%%%%%%%%%%%%%%%%%%%%%%%%%%%%%%%%%
%%%%%%%%%%%%%%%%%%%%%%%%%%
\section{Introduction}
In 1954, L. Carlitz [4, Theorem 5, p.144] gave a quite formula to
compute the number of points in the intersection of two quadrics.
His formula depended mainly on the capacity to write the first
quadric as a sum of non-degenerate quadrics $\mathcal{Q}_i$, the
number of indeterminates of each quadric $\mathcal{Q}_i$, and the
number of quadrics $\mathcal{Q}_i$ with odd indeterminates. Secondly
on the capacity to write the second quadrics as linear combination
of the
quadrics $\mathcal{Q}_i$.\\
In the same year A. Weil [17, p.348]  extend the result of L.
Carlitz, to two arbitrary quadrics. His formula depended on the
ability to diagonalize the forms in the pencil defined by the two
quadrics. We think that these sophisticated methods are concretly
difficult to apply, in order to find an interesting upper bound for
the number of
points in the intersection of two quadrics.\\
In 1975, W. M. Schmidt [14., Lemma 3C, p. 175], was the first the
give an explicit upper bound for the number of intersection of two
hypersurfaces by using some results on the theory of resultants.
Therefore for two quadrics, his estimate leads to the following
result:
$$| \mathcal{Q}_1\cap \mathcal{Q}_2|\le 2
(4q^{n-2}+4\pi_{n-3})+\frac{7}{q-1}.$$ But he was convinced that its
estimate was not the best possible.\\ In 1986, A. A. Bruen and J. W.
P. Hirschfeld [3, pp. 218-220] computed the number of points in the
intersection of two quadrics only for non-degenerate quadrics in a
projective space of odd dimension. Even in these cases, theirs
results don't satisfy totally our hope.\\
In 1992, Y. Aubry [1, Theorem 3, p. 11] inspired by the techniques
of W. M. Schmidt and the upper bound for hypersurfaces of
Tsfasman-Serre-S\o rensen [15], [16, chap.2,pp.7-10], improved the
above upper bound. We prefer to write his result in this way:
$$| \mathcal{Q}_1\cap \mathcal{Q}_2|\le 2
(4q^{n-2}+\pi_{n-3})+\frac{1}{q-1}.$$ In 1999, D. B. Leep and L. M.
Schueller [13, p.172] improved  the result of Y. Aubry under the
condition that the pair of quadrics have full order. This is a too
much restrictive condition. These are their bounds: $|
\mathcal{Q}_1\cap \mathcal{Q}_2|\le LS(n, q)$ with

\begin{equation*}
       LS(n, q) =  \begin{cases}

 2q^{n-2}+\pi_{n-3}+2q^{\frac{n-1}{2}}-q^{\frac{n-3}{2}} &  \text{if}\quad n+1 \geq 4\qquad \text{and}\qquad \text{even} \\

 2q^{n-2}+\pi_{n-3}+q^{\frac{n}{2}}& \text{if}\quad n+1 \geq 5\quad \text{and}\qquad
  \text{odd}.
            \end{cases}
\end{equation*}
In 2008, F. A. B. Edoukou [6] gave the best upper bound for the
number of points in the intersection of two quadrics with no common
plane in ${\mathbb{P}}^3(\mathbb{F}_q)$, when one of them is
non-degenerate (hyperbolic or elliptic) or degenerate of rank 3. He
deduced that
$4q+1$ is  the optimal bound.\\
In 2009, F. A. B. Edoukou [7] study the same problem in
${\mathbb{P}}^4(\mathbb{F}_q)$ when one of the two quadrics is
non-degenerate quadric (i.e. parabolic), degenerate of rank 4 or 3.
He deduced that $4q^2+q+1$ is the best upper bound and he formulated
the following conjecture:\\\\
\textbf{Conjecture [7, p.31]:} Let $\mathcal{Q}_1$ and
$\mathcal{Q}_2$ two quadrics in ${\mathbb{P}}^n(\mathbb{F}_q)$ with
no common hyperplane. Then
$$| \mathcal{Q}_1\cap \mathcal{Q}_2|\le 4q^{n-2}+\pi_{n-3}$$
And this bound is the best possible.\\\\
In this paper we will give a proof of this conjecture. The paper has
been organized as follows. First of all we recall some generaties on
quadrics. Secondly, based on what has been done for the study of
intersection of quadrics in the projective of three and four
dimension, and section of two quadrics with full order in higher
dimension, we prove the conjecture. Finally we will formulate
another conjecture for the intersection of quadrics, and a
generalization of the above conjecture for a projective algebraic
set of dimension $s$ and degree $d$ in a conjectural way:
$$|X(\mathbb{F}_q)|\le dq^s+\pi_{s-1}.$$ We will discuss on this
conjecture.
%%%%%%%%%%%%%%%%%%%%%%%%%%%%%%%%%%%%%%%%%%%%%%%%%%%%%%%%%%%%%%%%%%%%%%%%%
%%%%%%%%%%%%%%%%%%%%%%%%%%%%%%%%%%%%%%%%%%%%%%%%%%%%%%%%%%%%%%%%%%%%%%%%%%%
%%%%%%%%%%%%%%%%%%%%%%%%%%%%%%%%%%%%%%%%%%%%%%%%%%%%%%%%%%%%%%%%%%%%%%%%%%
%%%%%%%%%%%%%%%%%%%%%%%%%%%%%%%%%%%%%%%%%%%%%%%%%%%%%%%%%%%%%%%%%%%%%%%%%%
%%%%%%%%%%%%%%%%%%%%%%%%%%%%%%%%%%%%%%%%%%%%%%%%%%%%%%%%%%%%%%%%%%%%%%%%%%
 \section{Generalities}
We denote by $\mathbb{F}_q$ the field with $q$ elements and
${\mathbb{P}}^{n}(\mathbb{F}_q)=\Pi_n$ the projective space of
n-dimension over the field $\mathbb{F}_q$. Then
$$\pi_n=\#{\mathbb{P}^{n}(\mathbb{F}_q)}=q^n+q^{n-1}+...+1.$$
Let $\mathcal{Q}$ be a quadric in $\mathbb{P}^n(\mathbb{F}_q)$(i.e.
$\mathcal{Q}=Z(F)$ where $F$ is a form of degree 2). The rank of
$\mathcal{Q}$, denoted r$(\mathcal{Q})$, is the smallest number of
indeterminates appearing in $F$ under any change of coordinate
system. The quadric $\mathcal{Q}$ is said to be degenerate if
r$(\mathcal{Q})<n+1$; otherwise it is non-degenerate. For
$\mathcal{Q}$ a degenerate quadric and r$(\mathcal{Q})$=r,
$\mathcal{Q}$ is a cone $\Pi_{n-r}\mathcal{Q}_{r-1}$ with vertex
$\Pi_{n-r}$ (the set of the singular points of $\mathcal{Q}$) and
base $\mathcal{Q}_{r-1}$ in a subspace $\Pi_{r-1}$ skew to
$\Pi_{n-r}$.
%%%%%%%%%%%%%%%%%%%%%%%%%%%%%%%%%%%%%%%
\begin{definition}\label{quadrique type} For $\mathcal{Q}=\Pi_{n-r}\mathcal{Q}_{r-1}$
a degenerate quadric with r$(\mathcal{Q})=r$,  $\mathcal{Q}_{r-1}$ is called the non-degenerate
 quadric associated to $\mathcal{Q}$.
The degenerate quadric $\mathcal{Q}$ will be said to be of
hyperbolic type (resp. elliptic, parabolic) if its associated
non-degenerate quadric is of that type.
\end{definition}
%%%%%%%%%%%%%%%%%%%%%%%%%%%%%%%%%%%%%%%%
\begin{definition}\lbrack 13, p.158\rbrack\
Let $\mathcal{Q}_{1}$ and $\mathcal{Q}_{2}$ be two quadrics. The
order $w(\mathcal{Q}_{1}, \mathcal{Q}_{2})$ of the pair
$(\mathcal{Q}_{1}$, $\mathcal{Q}_{2})$, is the minimum number of
variables, after invertible linear change of variables, necessary to
write $\mathcal{Q}_{1}$ and $\mathcal{Q}_{2}$. When
$w(\mathcal{Q}_{1}, \mathcal{Q}_{2})=t+1$, we assume that
$\mathcal{Q}_1$ and $\mathcal{Q}_2$ are defined with the first $t+1$
indeterminates $x_0, \ldots, x_t$ and we define
$\mathbb{E}_{t}(\mathbb{F}_q)=\{x\in \mathbb{P}^n(\mathbb{F}_q)\
\vert x_{t+1}=...=x_n=0\}$.
\end{definition}
%%%%%%%%%%%%%%%%%%%%%%%%%%%%%%%%%%%%%%%%%
%\begin{definition} For any projective algebraic variety
%$\mathcal{V}$, the maximum dimension g$(\mathcal{V}$) of linear
%subspaces lying on $\mathcal{V}$, is called the projective index of
%$\mathcal{V}$. The largest dimensional spaces contained in
%$\mathcal{V}$ are called the generators of $\mathcal{V}$.
%\end{definition}
\begin{lemma}\label{Serre quadrique}\lbrack 5, pp.70-71 \rbrack  \
Let $\mathcal{Q}_1$ and $\mathcal{Q}_2$ be two quadrics in $\mathbb{P}^n(\mathbb{F}_q)$
and $l$ an integer such that $1\le l\le n-1$. %\\
Suppose that $w(\mathcal{Q}_1,\mathcal{Q}_2)=n-l+1$ (i.e. there exists a linear transformation
such that $\mathcal{Q}_1$ and $\mathcal{Q}_2$ are defined with the indeterminates
 $x_0, x_1,...,x_{n-l}$) and $\vert \mathcal{Q}_1\cap \mathcal{Q}_2\cap{\mathbb{E}}_{n-l}(\mathbb{F}_q)\vert \le m$.
 %where $\mathbb{E}_{n-l}(\mathbb{F}_q)=\{x\in \mathbb{P}^n(\mathbb{F}_q)\ \vert x_{n-l+1}=...=x_n=0\}$.
Then $$\vert \mathcal{Q}_1\cap \mathcal{Q}_2\vert \le
mq^{l}+\pi_{l-1}.$$ This bound is optimal as soon as $m$ is optimal
in $\mathbb{E}_{n-l}(\mathbb{F}_q)$.
\end{lemma}
%%%%%%%%%%%%%%%%%%%%%%%%%%%%%%%%%%%%%%%%%%%%%%%%%%%%%%%%%%%%%%%%%%%%%%%%%%%%%%%%%%%%%%%%%%%
%%%%%%%%%%%%%%%%%%%%%%%%%%%%%%%%%%%%%%%%%%%%%%%%%%%%%%%%%%%%%%%%%%%%%%%%%%%%%%%%%%%%%%%%%
%%%%%%%%%%%%%%%%%%%%%%%%%%%%%%%%%%%%%%%%%%%%%%%%%%%%%%%%%%%%%%%%%%%%%%%%%%%%%%%%%%%%%%%%%%%
%%%%%%%%%%%%%%%%%%%%%%%%%%%%%%%%%%%%%%%%%%%%%%%%%%%%%%%%%%%%%%%%%%%%%%%%%%%%%%%%%%%%%%%%%%%
%%%%%%%%%%%%%%%%%%%%%%%%%%%%%%%%%%%%%%%%%%%%%%%%%%%%%%%%%%%%%%%%%%%%%%%%%%%%%%%%%%%%%%%%%%%
\section{Resolution of the conjecture}
Let us state the following result which will be very useful in the
proof of the below theorem.
\begin{lemma}\label{cardinal hermitienne}
For $n\ge 3$,  $LS(n, q)\le 4q^{n-2}+\pi_{n-3}$.
\end{lemma}
\textbf{Proof:} An easy computation.
\begin{theorem}\label{BSupquadrique}Let $\mathcal{Q}_1$ and $\mathcal{Q}_2$ be two quadrics in
${\mathbb{P}}^n(\mathbb{F}_q)$ with no common hyperplane. Then
$$| \mathcal{Q}_1\cap \mathcal{Q}_2|\le 4q^{n-2}+\pi_{n-3}$$
And this bound is the best possible.
\end{theorem}
%%%%%%%%%%%%%%%%%%%%%%%%%%%%%%%%%%%%%%%%%%%%%%%%%%%%%%%%%%%%%%%%%%%%%%%%%%%%%%%%%%%%%%%%%%%%%%%%
%%%%%%%%%%%%%%%%%%%%%%%%%%%%%%%%%%%%%%%%%%%%%%%%%%%%%%%%%%%%%%%%%%%%%%%%%%%%%%%%%%%%%%%%%%%%%%%%
%%%%%%%%%%%%%%%%%%%%%%%%%%%%%%%%%%%%%%%%%%%%%%%%%%%%%%%%%%%%%%%%%%%%%%%%%%%%%%%%%%%%%%%%%%
\textbf{Proof:} It has been proved in [6, 7] that the above
conjecture is true for $n=3$ and 4. We will make an induction
reasoning on the dimension of the projective space.\\
Let us suppose that the conjecture is true for $i$, $4\le i \le
n-1$; that means for two quadrics $\mathcal{Q}_1$ and
$\mathcal{Q}_2$ in ${\mathbb{P}}^i(\mathbb{F}_q)$ with no common
hyperplane we have
$$\mathrm{for}\qquad 4\le
i\le n-1,\qquad | \mathcal{Q}_1\cap \mathcal{Q}_2|\le
4q^{i-2}+\pi_{i-3}.$$
Suppose now that $\mathcal{Q}_1$ and
$\mathcal{Q}_2$ are two quadrics in ${\mathbb{P}}^n(\mathbb{F}_q)$
with no common hyperplane. Let us consider the order
$w(\mathcal{Q}_{1}, \mathcal{Q}_{2})$ of the two quadrics. We have
either $w(\mathcal{Q}_{1}, \mathcal{Q}_{2})=n+1$
or $w(\mathcal{Q}_{1}, \mathcal{Q}_{2})\le n$.\\
If $w(\mathcal{Q}_{1}, \mathcal{Q}_{2})=n+1$, from the result of D.
B. Leep and L. M. Schueller, we deduce that $| \mathcal{Q}_1\cap
\mathcal{Q}_2|\le LS(n, q)$. Therefore from the above Lemma we
deduce
that the conjecture is true.\\
If $w(\mathcal{Q}_{1}, \mathcal{Q}_{2})\le n$, then there is an
integer $l$, $1\le l \le n-1$ such that $w(\mathcal{Q}_{1},
\mathcal{Q}_{2})=n-l+1$. Therefore there exists a linear
transformation such that $\mathcal{Q}_1$ and $\mathcal{Q}_2$ are
defined with the indeterminates $x_0, x_1,...,x_{n-l}$. And from the
above hypothese we deduce that: $$\vert \mathcal{Q}_1\cap
\mathcal{Q}_2\cap{\mathbb{E}}_{n-l}(\mathbb{F}_q)\vert \le
4q^{n-l-2}+\pi_{n-l-3}.$$ Then by applying Lemma \ref{Serre
quadrique} with $m=4q^{n-l-2}+\pi_{n-l-3}$, we deduce that $$|
\mathcal{Q}_1\cap \mathcal{Q}_2|\le 4q^{n-2}+\pi_{n-3}.$$ Let now $
\mathcal{Q}_1$ and $\mathcal{Q}_2$ be two quadrics of
${\mathbb{P}}^n(\mathbb{F}_q)$ defined by the equations $f_1(x_0,
x_1,...,x_n)=x_0^2+x_1^2-x_2^2$, and $f_2(x_0, x_1,...,x_n)=x_0x_1$.
In [Edoukou IEEE], it has been proved that (with $q$ odd), $|
\mathcal{Q}_1\cap \mathcal{Q}_2\cap {\mathbb{E}}_3(\mathbb{F}_q)|=
4q+1$. Therefore by using the above Lemma, we deduce that, we have
$| \mathcal{Q}_1\cap \mathcal{Q}_2|= 4q^{n-2}+\pi_{n-3}$. If $
\mathcal{Q}_1$ and $\mathcal{Q}_2$ are defined the equations
$f_1(x_0, x_1,...,x_n)=(x_0+x_1)x_2+x_2^2$, and $f_2(x_0,
x_1,...,x_n)=(x_2+x_0)x_1+x_1^2$, we have $| \mathcal{Q}_1\cap
\mathcal{Q}_2\cap {\mathbb{E}}_2(\mathbb{F}_q)|= 4$. And by the
preceding reasoning we have $| \mathcal{Q}_1\cap \mathcal{Q}_2|=
4q^{n-2}+\pi_{n-3}$.

\begin{Remark}
It has been proved in [8] that when one of the two quadrics
$\mathcal{Q}_1$ or $\mathcal{Q}_2$ is non-degenerate (hyperbolic,
elliptic, parabolic), we have $$| \mathcal{Q}_1\cap
\mathcal{Q}_2|\le EH(n, q)$$ with
\begin{equation*}
       EH(n, q) =  \begin{cases}

 2q^{n-2}+\pi_{n-3}+2q^{\frac{n-1}{2}}-q^{\frac{n-3}{2}} &  \text{if}\quad n+1 \geq 4\quad \text{and}\qquad \text{even} \\

 2q^{n-2}+\pi_{n-3}+q^{\frac{n-3}{2}} &  \text{if}\quad n+1 \geq 4\quad \text{and}\qquad \text{even} \\

 2q^{n-2}+\pi_{n-3}+2q^{\frac{n-2}{2}}& \text{if}\quad n+1 \geq 5\quad \text{and}\qquad
  \text{odd}.
            \end{cases}
\end{equation*}
These bounds are the best possible. Therefore, in this particular
case, these bounds are better than the one of D. B. Leep and L. M.
Schueller [13, p.172].
\end{Remark}
\section{Conjectures on the number of points of algebraic sets}
Here we will establish two conjectures concerning the number of
points in the intersection of two degenerate quadrics with no common
hyperplane and secondly for the number of points of a general algebraic set.\\
If $\mathcal{X}=\Pi_{n-1}\mathcal{P}_{0}$ is the degenerate quadric
of rank r=1 (i.e. a repeated hyperplane) and $\mathcal{Q}$ any
quadric not containing $\mathcal{X}$ we can easily prove that $|
\mathcal{X}\cap \mathcal{Q}|\le 2q^{n-2}+\pi_{n-3}$. And this upper
bound is optimal if we take for example $ \mathcal{X}$ and
$\mathcal{Q}$ defined by the equations $f_{\mathcal{X}}(x_0,
x_1,...,x_n)={x_2}^2$, and $f_{\mathcal{Q}}(x_0,
x_1,...,x_n)=x_0x_1+{x_2}^2$.\\
If $\mathcal{X}=\Pi_{n-2}\mathcal{E}_{1}$ is the degenerate quadric
of rank $r=2$ of elliptic type, for any quadric $\mathcal{Q}$ we
have $|\mathcal{X}\cap \mathcal{Q}|\le \pi_{n-2}$. And this upper
bound is reached if we take for example $\mathcal{X}$ and
$\mathcal{Q}$ defined by the equations $f_{ \mathcal{X}}(x_0,
x_1,...,x_n)=f(x_0,x_1)$ ($f$ is irreducible), and $f_{
\mathcal{Q}}(x_0,
x_1,...,x_n)=x_0x_1$.\\
We also know from the examples in the previous section that when
$\mathcal{X}=\Pi_{n-r}\mathcal{Q}_{r-1}$ is a degenerate quadric of
rank $r=2$ of hyperbolic type or rank $r=3$, we can find another
degenerate quadric such that $|
\mathcal{X}\cap \mathcal{Q}|=4q^{n-2}+\pi_{n-3}$.\\
Let us now suppose that $r\ge 4 $. We have the following conjecture
\begin{Conjecture}
Let $\mathcal{X}=\Pi_{n-r}\mathcal{Q}_{r-1}$ be a degenerate quadric
of rank r$(\mathcal{X})=r$ and $\mathcal{Q}$ another degenerate
quadric.  If $r\ge 4$, then
$$|\Pi_{n-r}\mathcal{Q}_{r-1} \cap
\mathcal{Q}|\le EH(r-1, q)q^{n-r+1}+\pi_{n-r}.$$ And this bound is
the best possible
\end{Conjecture}
Let us remark that if the conjecture is true, that means for every
functional code of order 2 defined on a degenerate quadric
$\mathcal{X}=\Pi_{n-r}\mathcal{Q}_{r-1}$, we know exactly its
minimum distance.\\\\
Based on the authors experiences on the study of the number of
points in the intersection of hypersufaces, we would like  to
propose the following conjecture.
\begin{Conjecture}\label{EdoukConj}
Let $X\subset \mathbb{P}^{n}(\mathbb{F}_q)$ be a projective
algebraic set of degree $d$ and dimension $s$. Then,
\begin{equation} \label{EdoukConj1}
|X(\mathbb{F}_q)|\le dq^s+\pi_{s-1}.
\end{equation}
\end{Conjecture}
\begin{Remark}
Let us remark that this conjecture is true if X is of codimension 1.
This is the Tsfasman-Serre-S\o rensen's upper bound for
hypersurfaces [15], [16, chap.2, pp.7-10].\\ The Theorem
\ref{BSupquadrique} is a particular case of this conjecture. In
fact, the intersection of the two quadrics with no common hyperplane
is an algebraic set of dimension $n-2$ and degree $4$.
\end{Remark}
\begin{Remark}There is also an upper bound concerning the number of points
of a projective algebraic set which was already known by G. Lachaud
since 1993 and proved by Lachaud[12, p.80], Boguslavsky[2, pp.288,
293-294] and Gorpade-Lachaud [9,pp.627-630]:
 $$|X(\mathbb{F}_q)|\le d\pi_{s}.$$
 In 1995, G. Lachaud [2, p. 292], [9, p.629] proposed the
 following conjecture on the upper bound for the number
 of points of a projective algebraic set:
\begin{equation} \label{LachConj2}
|X(\mathbb{F}_q)|\le d(\pi_s-\pi_{2s-n})+\pi_{2s-n}.
\end{equation}
Our bound (\ref{EdoukConj1}) is more better than the one of G.
Lachaud (\ref{LachConj2}). In [5, 6, 7, 8] we can find several cases
where the Conjecture \ref{EdoukConj} is true and where the upper
bound of G. Lachaud (\ref{LachConj2}) seems to be too large.
\end{Remark}
%%%%%%%%%%%%%%%%%%%%%%%%%%%%%%%%%%%%%%%%%%%%%%%%%%%%%%%%%%%%%%%%%%%%%%%%%%%%%%%%%%%%%%
%%%%%%%%%%%%%%%%%%%%%%%%%%%%%%%%%%%%%%%%%%%%%%%%%%%%%%%%%%%%%%%%%%%%%%%%%%%%%%%%%%%%%%
%%%%%%%%%%%%%%%%%%%%%%%%%%%%%%%%%%%%%%%%%%%%%%%%%%%%%%%%%%%%%%%%%%%%%%%%%%%%%%%%%%%%%%
%%%%%%%%%%%%%%%%%%%%%%%%%%%%%%%%%%%%%%%%%%%%%%%%%%%%%%%%%%%%%%%%%%%%%%%%%%%%%%%%%%%%%%
\textbf{References}\\
{\footnotesize \lbrack 1\rbrack  \ Y. Aubry, Reed-Muller codes
associated to projective varieties, in "Coding Theory and Algebraic
Geometry, Luminy, 1991," Lecture Notes in Math. Vol. 1518
Springer-Verlag, Berlin, 4-17, 1992.\\
\lbrack 2\rbrack  \ M. Boguslavsky, On the number of points of
polynomial systems. Finite Fields and Theirs Applications 3, (1999),
287-299.\\
\lbrack 3\rbrack \ A. A. Bruen and J. W. P. Hirschfeld,
Intersections in projective space I: Combinatorics. Mathematische
Zeitschrift, Vol. 193, (1986), 215-226.\\
\lbrack 4\rbrack  \ L. Carlitz, Pairs of quadratic equations in a
finite field. America Journal of Mathematics, Vol. 76, 137-154,
(1954).\\
 \lbrack 5 \rbrack  \ F. A. B. Edoukou,
Codes correcteurs d'erreurs construits \`a partir des vari\'et\'es
alg\'ebriques. Ph.D Thesis, Universit\'e de la M\'editerran\'ee (Aix-Marseille II), France, (2007).\\
\lbrack 6\rbrack  \ F. A. B. Edoukou, Codes defined by forms of
degree 2 on quadric surfaces. IEEE Trans. Inf. Vol 54 N0.(2), 860-864, (2008).\\
\lbrack 7\rbrack \ F. A. B. Edoukou, Codes defined by forms of
degree 2 on quadric varieties in $\mathbb{P}^{4}(\mathbb{F}_q)$.
Contemporary Mathematics, Proceedings of the 11 th
Conference AGC2T, Marseille, France, November 2007,  21-32, (2009).\\
\lbrack 8\rbrack  \ F. A. B. Edoukou, A. Hallez, F. Rodier and L.
Storme, A study of intersections of quadrics having applications on
the small weight codewords of the functional codes
$C_2(\mathcal{Q})$, $\mathcal{Q}$ a non-singular quadric. 19 pages,
 Submitted to Journal of Pure and Applied Algebra (2009).\\
\lbrack 9\rbrack \ S. R. Ghorpade and G. Lachaud, Etale cohomology,
Lefstchetz theorems and number of points of singular varieties over
finite fields. Moscow Mathematical Journal, Vol. 2, N0.3, (2002),
589-631.\\
\lbrack 10\rbrack  \ J. W. P. Hirschfeld, Projective Geometries Over
Finite Fields (Second Edition)
 Clarendon  Press. Oxford 1998.\\
\lbrack 11\rbrack  \ J. W. P. Hirschfeld,  General Galois Geometries, Clarendon press. Oxford 1991.\\
\lbrack 12\rbrack \ G. Lachaud, Number of points of plane sections
and linear codes defined on
 algebraic varieties;  in " Arithmetic, Geometry, and Coding Theory ". (Luminy, France, 1993),
 Walter de Gruyter, Berlin-New York, (1996), 77-104. \\
\lbrack 13\rbrack \ D. B. Leep and L. M. Schueller, Zeros of a pair
of quadric forms defined over
  finite field. Finite Fields and Their Applications 5 (1999), 157-176.\\
\lbrack 14\rbrack \ W. M. Schmidt, Equations over finite fields. An
Elementary  Approach, Lecture Notes in Maths 536 (1975).\\
\lbrack 15\rbrack \ J. -P. Serre, Lettre \`a M. Tsfasman, In
"Journ\'ees Arithm\'etiques de Luminy (1989)", Ast\'erisque
198-199-200 (1991), 3511-353.\\
\lbrack 16\rbrack  \ A. B. S\o rensen, Rational points on
hypersurfaces, Reed-Muller codes and
 algebraic-geometric codes. Ph. D. Thesis, Aarhus, Denmark, 1991.\\
\lbrack 17\rbrack \ A. Weil, Footnote to a recent paper. American
Journal of Mathematics, Vol. 76, 347-350, (1954).}
\end{document}